\documentclass[a4paper,fleqn]{cas-sc}

\usepackage[numbers,sort]{natbib}

\usepackage{subcaption} %
\usepackage{siunitx} %
\usepackage{lipsum} %
\usepackage[mathlines]{lineno} %
\usepackage[ruled,vlined]{algorithm2e} %
\usepackage[nonumberlist,acronym]{glossaries} %
\usepackage{mathtools}

\graphicspath{{figs/}{figs/soi_dotted/}}

\DeclareSIUnit{\AU}{AU}
\DeclareSIUnit{\day}{d}
\DeclareSIUnit{\deg}{deg}
\DeclareSIUnit{\year}{yr}

\newcommand{\di}{{\rm d}}

\newcommand{\Eq}[1]{Eq.~\eqref{#1}}

\newcommand{\Fig}[1]{Fig.~\ref{#1}}

\newcommand{\Figsand}[2]{Figs.~\ref{#1} and~\ref{#2}}
\newcommand{\Tab}[1]{Table~\ref{#1}}

\newcommand{\Sec}[1]{Section~\ref{#1}}

\newcommand{\tr}{\textcolor{red}}

\newcommand{\ie}{i.\,e.,~}

\newcommand{\eg}{e.\,g.,~}

\newacronym{ABM}{ABM}{Adams-Bashforth-Moulton}
\newacronym{ADS}{ADS}{automatic domain splitting}
\newacronym{BC}{BC}{ballistic capture}
\newacronym{CR3BP}{CR3BP}{circular restricted three-body problem}
\newacronym{DAER}{DAER}{Department of Aerospace Science and Technology}
\newacronym{DOPRI8}{DOPRI8}{Dormand-Prince 8th-order embedded Runge-Kutta method}
\newacronym{EOM}{EoM}{equations of motion}
\newacronym{ER3BP}{ER3BP}{elliptic restricted three\tr{-}body problem}
\newacronym{ERC}{ERC}{European Research Council}
\newacronym{IC}{IC}{initial condition}
\newacronym{JPL}{JPL}{Jet Propulsion Laboratory}
\newacronym{LD}{LD}{Lagrangian descriptor}
\newacronym{LCS}{LCS}{Lagrangian coherent structure}
\newacronym{NAIF}{NAIF}{Navigation and Ancillary Information Facility}
\newacronym{NASA}{NASA}{National Aeronautics and Space Administration}
\newacronym{PECE}{PECE}{predictor-corrector}
\newacronym{RAAN}{RAAN}{right ascension of the ascending node}
\newacronym{RHS}{RHS}{right-hand side}
\newacronym{RK}{RK}{Runge-Kutta}
\newacronym{RSS}{RSS}{root sum square}
\newacronym{SOI}{SOI}{sphere of influence}
\newacronym{SRP}{SRP}{solar radiation pressure}
\newacronym{UOM}{uom}{unit of measurement}
\newacronym{VSVO}{VSVO}{variable-step variable-order}
\newacronym{WSB}{WSB}{weak stability boundary}

\def\tsc#1{\csdef{#1}{\textsc{\lowercase{#1}}\xspace}}
\tsc{WGM}
\tsc{QE}

\begin{document}
\let\WriteBookmarks\relax
\def\floatpagepagefraction{1}
\def\textpagefraction{.001}

\shorttitle{Qualitative study of ballistic capture at Mars via Lagrangian descriptors}

\shortauthors{A. Quinci et al.}  

\title [mode = title]{Qualitative study of ballistic capture at Mars via Lagrangian descriptors}

\author[1]{Alessio Quinci}[type=editor,
                      orcid=0000-0002-7792-4581]

\cormark[1]
\ead{alessio.quinci@mail.polimi.it}
\credit{Investigation, Methodology, Writing -- Original draft preparation}
\affiliation[1]{organization={Department of Aerospace Science and Technology, Politecnico di Milano},
                addressline={Via La Masa 34}, 
                city={Milano},
                citysep={}, %
                postcode={20156}, 
                country={Italy}}

\cortext[1]{Corresponding author}

\author[1]{Gianmario Merisio}[orcid=0000-0001-8806-7952]
\cormark[2]
\ead{gianmario.merisio@polimi.it}
\credit{Methodology, Writing -- Original draft preparation}
\cortext[2]{Principal corresponding author}

\author[1]{Francesco Topputo}[orcid=0000-0002-5369-6887]
\ead{francesco.topputo@polimi.it}
\credit{Conceptualization of this study, Methodology, Draft revision, Funding acquisition}

\begin{abstract}
Lagrangian descriptors reveal the dynamical skeleton governing transport mechanisms of a generic flow. In doing so, they unveil geometrical structures in the phase space that separate regions with different qualitative behavior. This work investigates to what extent \acrlongpl*{LD} provide information about non-Keplerian motion in Mars proximity, which is modeled under the planar \acrlong*{ER3BP}. We propose a novel technique to reveal \acrlong*{BC} orbits extracting separatrices of the phase space highlighted by \acrlong*{LD} scalar fields. The Roberts' operator to approximate the gradient is used to detect the edges in the fields. Results demonstrate the chaos indicator ability to distinguish sets of \acrlongpl*{IC} exhibiting different dynamics, including \acrlong*{BC} ones. Separatrices are validated against reference \acrlong*{WSB} derived on similar integration intervals. Compared to other techniques, \acrlongpl*{LD} provide dynamics insight bypassing the propagation of the variational equations.
\end{abstract}

\begin{highlights}
\item Lagrangian descriptors highlight regions with different dynamical behavior.
\item Geometrical singularities are extracted with Roberts’ edge detection method.
\item Dynamics separatrices are inspected against the weak stability boundary.
\item Insight about the ballistic capture problem in the proximity of Mars is provided.
\item Lagrangian descriptors are convenient for designing ballistic capture orbits.
\end{highlights}

\begin{keywords}
	Elliptic restricted 3-body problem \sep Lagrangian descriptors \sep Ballistic capture at Mars \sep Edge detection 
\end{keywords}

\maketitle

\section{Introduction} \label{sec:introduction}
\Gls*{BC} orbits are low-energy transfers that allow temporary capture about a planet exploiting the natural dynamics, thus without requiring maneuvers \cite{topputo2015earth}. \tr{Compared to Keplerian solutions, they are cheaper and more versatile from the operational perspective at the expense of longer transfer times.} \gls*{BC} orbits are bounded by the \gls*{WSB} \cite{belbruno1993sun,belbruno2000calculation,circi2001dynamics,topputo2015earth}. After being initially conceived as a fuzzy boundary region in the Sun--Earth--Moon system \citep{belbruno1987lunar,belbruno1990ballistic}, the \gls*{WSB} was algorithmically defined in \cite{belbruno2004capture}. The definition was later extended in \cite{garcia2007note,topputo2009computation,silva2012applicability}. A formal definition and a technique for its derivation were proposed in \cite{hyeraci2010method}.

\tr{Approaches currently known for designing \gls*{BC} orbits are: i) the technique stemmed from invariant manifolds \cite{conley1969ultimate,topputo2005low,belbruno2010weak}, ii) the method based on stable sets manipulation \cite{hyeraci2010method,topputo2009computation,luo2014constructing,luo2015analysis}, iii) the Hamiltonian approach taking advantage of canonical transformations \cite{carletta2019design}, and iv) the multiple shooting technique to solve a sequence of three-point boundary value problems \cite{merisio2023engineering}. The first methodology gives insights into the dynamics but it is only applicable to autonomous systems (\eg the \acrlong*{CR3BP}), while the others can be applied to more representative, non-autonomous models.} Lately, the variational theory for \glspl*{LCS} \cite{haller2011variational,haller2015lagrangian}, and the Taylor differential algebra \cite{wittig2015propagation} were applied to derive \gls*{BC} orbits and the \gls*{WSB} more efficiently \cite{manzi2021flow,caleb2022stable,tyler2022improved}. Alternatively, \glspl*{LD} can be exploited. They reveal separatrices, so providing a qualitative description of the dynamics and highlighting the geometrical template of phase space structures even for systems with generic time dependence \cite{jimenez2009distinguished,mancho2013lagrangian,lopesino2017theoretical,raffa2023finding,bernardini2022exploiting}.

The goal of the paper is to study to what extent \glspl*{LD} inform about the \gls*{BC} mechanism and aid in the design of \gls*{BC} orbits. We provide a characterization of the dynamics in the Mars proximity modeled under the planar \gls*{ER3BP}. The geometrical structures featured by \gls*{LD} scalar fields are extracted through an edge detection algorithm based on the Roberts' method \cite{davis1975survey}. \tr{Specifically, Roberts' operator is used to approximate the gradient of the field \cite{roberts1963machine}}. The separatrices are inspected against the \gls*{WSB} derived on similar integration intervals. For a coherent comparison, the particle stability definition is modified to relax the geometrical constraint on the number of completed revolutions \cite{hyeraci2010method,luo2014constructing}. Results show a strong correlation between extracted separatrices and the \gls*{WSB}, particularly when the geometrical structures governing the transport mechanisms emerge. Eventually, capture sets at Mars are identified in the intricate plot of separatrices.

The remainder of the paper is organized as follows. In \Sec{sec:eqofmotion}, the dynamical model is described. The methodology is discussed in \Sec{sec:methodology}. Results are shown in \Sec{sec:results}. Eventually, conclusions are drawn in \Sec{sec:conclusion}.

\section{Equations of motion} \label{sec:eqofmotion}
The planar \gls*{ER3BP} describes the motion of a massless particle moving under the gravitational attraction of two primary bodies $P_{1}$ (the Sun) and $P_{2}$ (Mars) without influencing their motion. The two primaries revolve on ellipses about their common barycenter, influenced only by their mutual attraction. The model is expressed in the synodic reference frame centered at the primaries barycenter. The synodic frame non-uniformly rotates and pulsates to keep their distance equal to one \cite{hyeraci2010method}. Let the mass parameter $\mu = m_{2}/(m_{1}+m_{2})$, where $m_{1}$ and $m_{2}$ are the masses of $P_{1}$ and $P_{2}$, respectively. The positions of $P_{1}$ and $P_{2}$ are (-$\mu$, 0) and (1-$\mu$, 0), respectively. The \gls*{EOM} are scaled such that the sum of $P_{1}$ and $P_{2}$ masses is set to one as well as their distance, and their period is scaled to $2\pi$ \cite{hyeraci2010method}. The true anomaly $f$ is designated as the independent variable of the system. The \gls*{EOM} read \cite{hyeraci2010method}
\begin{linenomath*}
\begin{align}
    \begin{split}
        {x}''-2{y}' &= \omega_{x} \\
        {y}''+2{x}' &= \omega_{y} 
    \end{split}
    \label{eq:er3bp}
\end{align}
\end{linenomath*}
where primes represent differentiation with respect to the true anomaly $ f $ that depends on the scaled time as \cite{hyeraci2010method}
\begin{linenomath*}
\begin{equation}
	\frac{\di f}{\di t} = \frac{(1+e_{p}\cos{f})^{2}}{(1-e_{p}^{2})^{3/2}}.
\end{equation}
\end{linenomath*}
In \Eq{eq:er3bp}, subscripts $(\cdot)_{x}$ and $(\cdot)_{y}$ denote the partial derivatives of the potential function $\omega$ defined as \cite{hyeraci2010method}
\begin{linenomath*}
\begin{equation}
    \omega \left(x,y,f\right)=\frac{1}{1+e_{p}\cos{f}} \left[ \frac{1}{2} \left(x^2+y^2\right) +\frac{1-\mu}{r_1} +\frac{\mu}{r_2} +\frac{1}{2}\mu(1-\mu) \right],
\end{equation}
\end{linenomath*}
with $r_{1} = \sqrt{(x+\mu)^2+y^2}$ and $ r_{2} = \sqrt{(x+\mu-1)^2+y^2}$ the distances of the particle from $P_1$ and $P_2$, respectively, while $e_{p}$ is the common eccentricity of the primaries. The Sun--Mars physical parameters used in this study are reported in \Tab{tab:parameters}. The \gls*{EOM} are integrated with a $8^{\rm th}$-order Runge--Kutta scheme with a $7^{\rm th}$-order embedded step-size control. The integration relative tolerance is set to $10^{-9}$ \cite{montenbruck2000satellite,verner2010numerically}.

\begin{table}[pos=tbp,width=0.70\textwidth,align=\centering]
    \caption{Sun--Mars physical parameters.} 
    \begin{tabular*}{0.70\textwidth}{@{}LLLLC@{}}
    \toprule
    \textbf{Parameter} & \textbf{Unit} & \textbf{Value} & \textbf{Description} & \textbf{Reference} \\
    \midrule
    $ \mu $ & - & \SI{3.226201e-7}{} & Mass parameter & \multirow{4}{*}{\cite{hyeraci2010method}} \\
    $ a_{p} $ & \si{\AU} & \SI{1.523688}{} & Primaries semi-major axis & \\ 
    $ e_{p} $ & - & \SI{0.093418}{} & Primaries eccentricity & \\
    $ R $ & \si{\kilo\meter} & \si{3397} & Mars mean equatorial radius & \\
    \cline{5-5}
    $ R_{\mathrm{SOI}} $ \rule{0pt}{10pt} & \si{\kilo\meter} & $ \SI{170}{} R $ & Mars \gls*{SOI} radius & \multirow{1}{*}{\cite{luo2014constructing}} \\
    \bottomrule
    \end{tabular*}
    \label{tab:parameters}
\end{table}

\section{Methodology} \label{sec:methodology}

{\color{red}
\subsection{Low-energy regime}
\gls*{BC} orbits can be identified as those solutions allowing for material transfer between interior and exterior realms \cite{conley1968low}. In the \gls*{CR3BP}, they could be quantitatively identified as trajectories with Jacobi constant just below that of the collinear Lagrangian points $L_1$ and $L_2$. In the Sun--Mars system under study, this means trajectories with $ C_J < C_{J1} = \num{3.000203} $ (\ie transfers between primary and secondary interior realms) and $ C_J < C_{J2} = \num{3.000202} $ (\ie transfers between interior and exterior realms) \cite{conley1968low}. On the other hand, when the Jacobi constant falls below that of the collinear Lagrangian point $L_3$ ($ C_J < C_{J3} = \num{3.000001} $), then high-energy transfers are found \cite{campagnola2010endgame,restrepo2017patched}. However, qualitative statement on allowed and forbidden regions are no longer possible in the elliptic problem since the Jacobi value becomes anomaly-dependent \cite{hyeraci2010method}.}

\subsection{Particle stability definition} \label{sec:stability-def}
Particle stability is inferred using an alternative formulation of the stable sets defined in \cite{luo2014constructing}. While propagating \glspl*{IC} in the non-dimensional, synodic reference frame, the particle non-dimensional distance $ r(f) $ and Kepler energy $ H(f) $ with respect to the target body $P_2$ are computed \cite{hyeraci2010method}. The following indications are used to classify stability: A) a particle escapes at $f=f_e$ if $ H(f_e)>0 \ \wedge \ r(f_e)>R_{\mathrm{SOI}} $; B) a particle impacts  the surface of the target at $f=f_i$ if $ r(f_i)<R $. Based on its dynamical behavior over the integration interval $ [f_{0}, f_{f}] $, a propagated trajectory is said to be: i) \textit{weakly stable} if the particle neither escape nor impact with the target, so belonging to the subset $ \mathcal{W}(f_{f}) $; ii) \textit{unstable} if the particle escapes from the target before $ f_{f} $, then condition A) is verified for $ f_{e} \in [f_{0}, f_{f}] $, so belonging to the subset $ \mathcal{X}(f_{f}) $; iii) \textit{crash} if the particle impacts with the target before $ f_{f} $, then condition B) is verified for $ f_{i} \in [f_{0}, f_{f}] $, so belonging to the subset $ \mathcal{K}(f_{f}) $. A capture set is defined as $ \mathcal{C}(f_{B},f_{F}) \coloneqq \mathcal{X}(f_{B}) \cap \mathcal{W}(f_{F}) $ where $ f_{B} < f_{0} $ (backward leg), and $ f_{F} > f_{0} $ (forward leg).

\tr{If neither the escape criteria nor the impact criteria are matched, then the orbit is considered weakly stable in the interval $ [f_0, f_f] $, independently of the type of orbit considered. The definition adopted in this study is the same used in \cite{hyeraci2010method,luo2014constructing}. The aforementioned alternative particle stability formulation only regards the count of revolutions about the target, which in this work is neglected.}

\subsection{Lagrangian descriptors}
By manipulating the definition given in \cite{mancho2013lagrangian}, we define the \gls*{LD} as
\begin{linenomath*}
\begin{equation}
	M(\mathbf{x}_0,f_0,f_{B},f_{F}) = \int_{f_{0}+f_{B}}^{f_{0}} |\mathcal{F}(\mathbf{x}(f))|^{\gamma} \di f + \int_{f_{0}}^{f_{0}+f_{F}} |\mathcal{F}(\mathbf{x}(f))|^{\gamma} \di f,
	\label{eq:LDdef}	
\end{equation}
\end{linenomath*}
where $ \mathbf{x} = [x,y,x',y'] $ is the state vector obtained by rearranging \Eq{eq:er3bp} as a four-dimensional, first-order system of ordinary differential equations $ \mathbf{x}' = \mathbf{f}(\mathbf{x},f) $. The integrand $ |\mathcal{F}(\mathbf{x}(f))|^\gamma$ in \Eq{eq:LDdef} is a bounded, positive quantity, while $\gamma$ is the exponent defining the norm \cite{mancho2013lagrangian}. In this study, we select $ \mathcal{F} \coloneqq \sqrt{(x')^2+(y')^2} $ and $\gamma = 1/2 $ because they highlight the geometrical structures of the phase space better than the other integrands and norms as in \cite{mancho2013lagrangian}. The \gls*{LD} field is then defined as $ \mathcal{M}(f_{0},f_{B},f_{F}) \coloneqq \{ M(\mathbf{x}_0,f_{0},f_{B},f_{F}) \mid \mathbf{x}_{0} \in \Omega\} $, where $ \Omega $ is the set containing the \glspl*{IC}. In practise, the \gls*{LD} is computed appending its integrand to the space state equations with a zero initial value, and propagating the extended dynamics. The integration of the extended dynamics is stopped at $ f_{i} $ if the particle impacts with the target body.

An abrupt change in the \gls*{LD} field yelds discontinuous derivatives along the direction transverse to the change. Such singularities coincide with phase space structures separating trajectories with different dynamics, so abrupt changes correspond to dynamics separatrices \cite{mancho2013lagrangian}. In \Eq{eq:LDdef}, $ M(\mathbf{x}_0,f_{0},0,f_{F}) $ isolates dynamics separatrices obtained propagating \glspl*{IC} forwards, thus they are linked to repelling \glspl*{LCS}. Conversely, $ M(\mathbf{x}_0,f_{0},f_{B},0)$ reveals separatrices backwards, so highlighting the attracting \glspl*{LCS} \cite{lopesino2017theoretical}.

\subsection{Extraction of separatrices}
The structures revealed by the \gls*{LD} field are extracted with an edge detection algorithm. Edge detection is an image processing technique usually exploited for finding boundaries of objects within images. An edge is defined as the locus of points where an abrupt change in intensity of the image occurs. Several edge detection algorithms are available (\eg Sobel, Prewitt, Roberts, Canny, and zero-cross methods) \cite{davis1975survey}. \tr{Roberts' operator appears to be the most effective in extracting edges from the \gls*{LD} field of the problem at hand \cite{roberts1963machine}.} 

Given the image of the scalar field, the algorithm finds edges at those points where the gradient magnitude of the image is larger than a sensitivity threshold $ \sigma $ provided as input. The gradient of the image is approximated by computing the sum of the squares of the differences between diagonal neighbors pixels \cite{davis1975survey,roberts1963machine}. The threshold value is tuned\footnote{Values too small ($ < 10^{-3}$) could generate false positives in the output binary image. The larger the final true anomaly, the larger is the threshold suggested to use. The trend is justified because changes in the \gls*{LD} value at the separatrices are stronger for longer propagations.} to show as many structures as possible associated to abrupt changes in the \gls*{LD} field. 

\subsection{Validation of separatrices}

\begin{figure}[tbp]
	\centering
	\includegraphics[width=\textwidth]{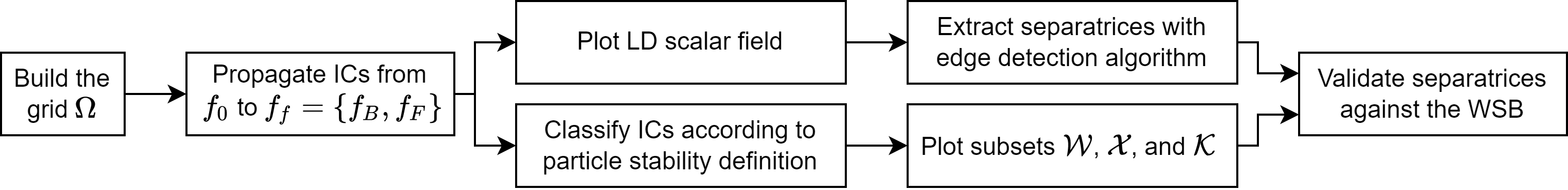}
	\caption{Validation workflow.}
	\label{fig:Method}
\end{figure}

The dynamics separatrices extracted from the \gls*{LD} field are expected to match the \gls*{WSB} computed on the same integration interval. The validation procedure devised to verify the correlation is outlined in \Fig{fig:Method}. Firstly, a uniform computational grid $ \Omega $ having $ 500 \times 500 $ points and centered at the target body is built over the square domain $ [-\varepsilon,\varepsilon] \times [-\varepsilon,\varepsilon] $, with $ \varepsilon = \num{6e-4} $. At $f_{0}$, the particle is assumed at the periapsis of an osculating prograde elliptic orbit about the target body with given eccentricity $ e_{0} = 0.9 $ (see \cite{hyeraci2010method} for more details). Secondly, \glspl*{IC} are propagated in the $ [f_{0}, f_{f}] $ interval. The \gls*{LD} values are computed and the \glspl*{IC} are allocated into the sets $\mathcal{W}$, $\mathcal{X}$, or $\mathcal{K}$ according to the stability definition discussed in \Sec{sec:stability-def}. Then, the separatrices are extracted with the edge detection algorithm. Finally, the patterns are inspected against the \gls*{WSB}.

\section{Results} \label{sec:results}

\begin{figure}[tbp]
	\centering
	\begin{subfigure}[t]{0.49\textwidth}
		\centering
		\includegraphics[width=\textwidth,clip,trim={0, 0, 0, 0}]{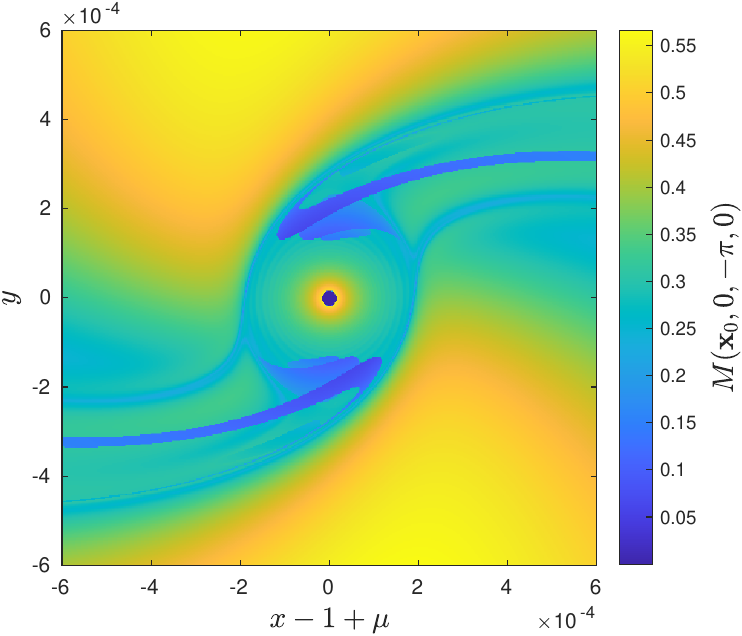}
		\caption{$\mathcal{M}(0,-\pi,0)$.}
		\label{fig:LD1}
	\end{subfigure}
	\begin{subfigure}[t]{0.49\textwidth}
		\centering
		\includegraphics[width=0.95\textwidth,clip,trim={1.45cm, 0, 1.45cm, 0}]{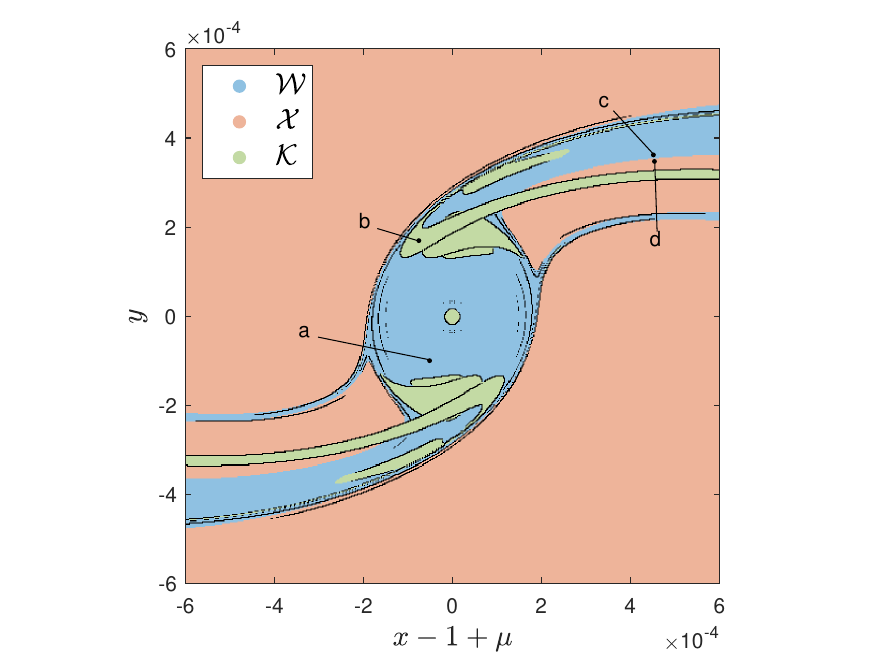}
		\caption{$\mathcal{M}(0,-\pi,0)$ separatrices and sets at $f_f=-\pi$; $ \sigma = \num{4e-3} $.}
		\label{fig:LD2}
	\end{subfigure}
	\begin{subfigure}[t]{0.49\textwidth}
		\centering
		\includegraphics[width=\textwidth,clip,trim={0, 0, 0, 0}]{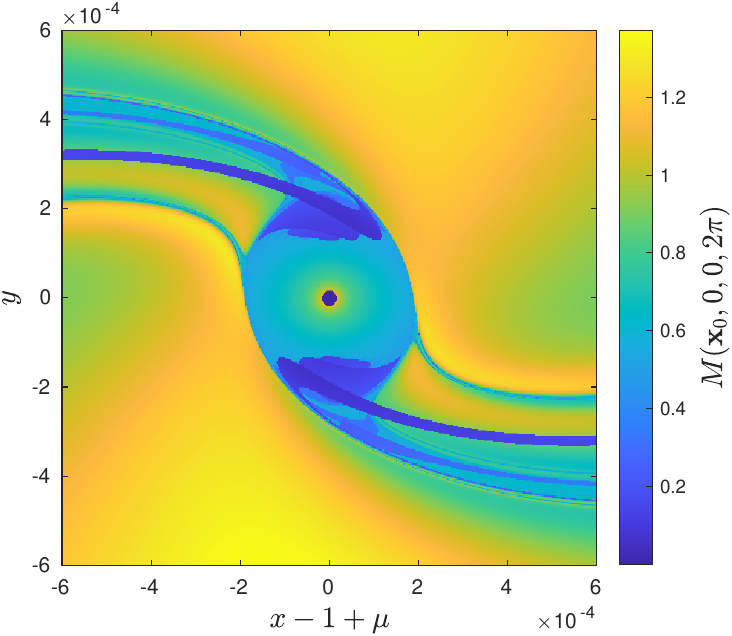}
		\caption{$\mathcal{M}(0,0,2\pi)$.}
		\label{fig:LD3}
	\end{subfigure}
	\begin{subfigure}[t]{0.49\textwidth}
		\centering
		\includegraphics[width=0.95\textwidth,clip,trim={1.45cm, 0, 1.45cm, 0}]{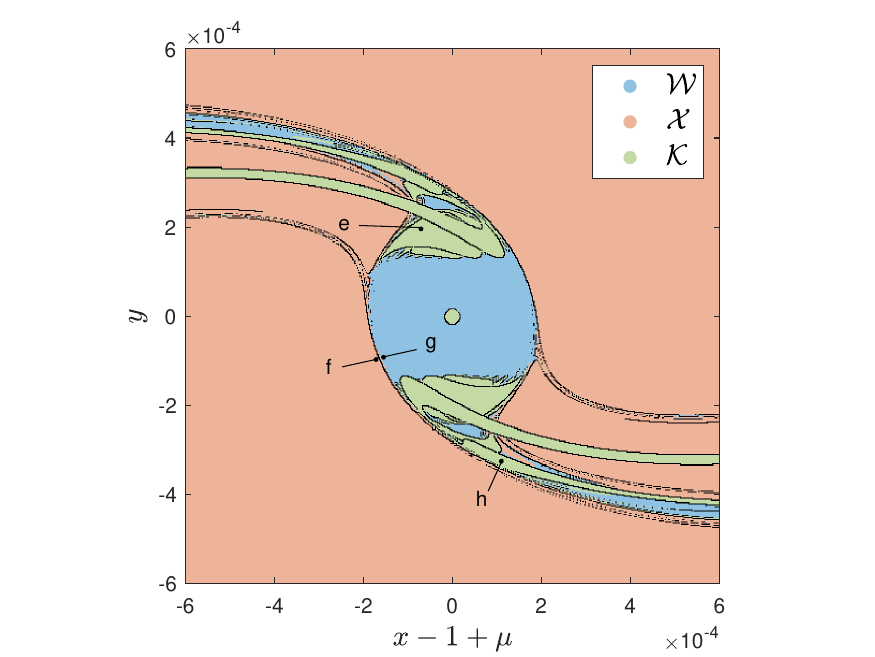}
		\caption{$\mathcal{M}(0,0,2\pi)$ separatrices and sets at $f_f=2\pi$; $ \sigma = \num{2e-2} $.}
		\label{fig:LD4}
	\end{subfigure}
	\caption{Validation of separatrices extracted from \gls*{LD} field through inspection against the \gls*{WSB}. Left: \gls*{LD} scalar field. Right: Extracted separatrices inspected against subsets $\mathcal{W}(f_f)$, $\mathcal{X}(f_f)$, and $\mathcal{K}(f_f)$.}
	\label{fig:LD}
\end{figure}

\begin{figure}[tbp]
	\centering
	\begin{subfigure}[t]{0.48\textwidth}
			\centering
			\includegraphics[width=\textwidth,clip,trim={2.1cm, 0, 2.1cm, 0}]{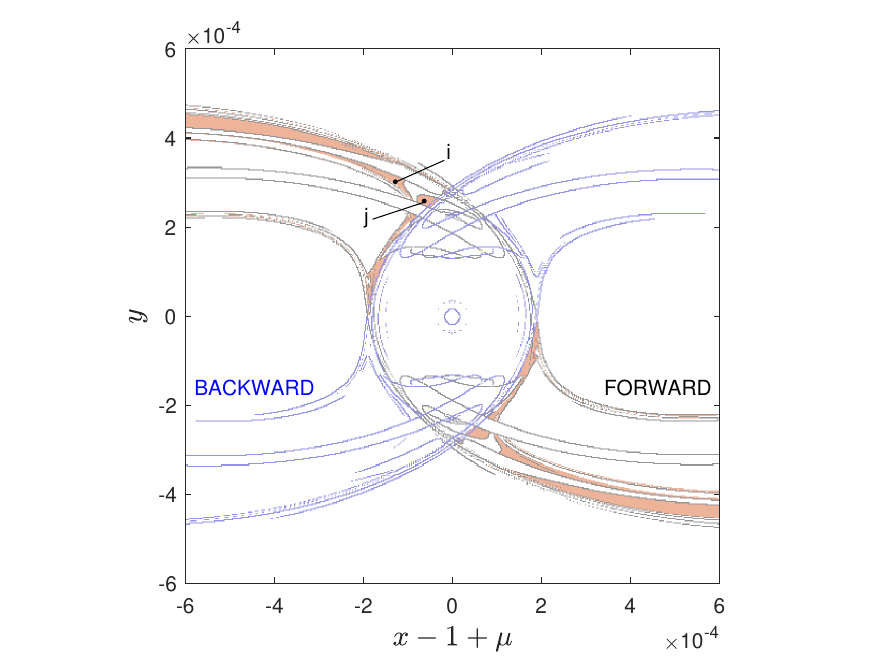}
			\caption{$\mathcal{M}(0,-\pi,3\pi/2)$ separatrices inspected against capture set $\mathcal{C}(-\pi,3\pi/2)$. Sensitivity threshold $ \sigma $ set to \num{4e-3} and \num{9e-3} for backward and forward propagations, respectively.}
			\label{fig:CaptureSets1}
	\end{subfigure} \quad
	\begin{subfigure}[t]{0.48\textwidth}
			\centering
			\includegraphics[width=\textwidth,clip,trim={2.1cm, 0, 2.1cm, 0}]{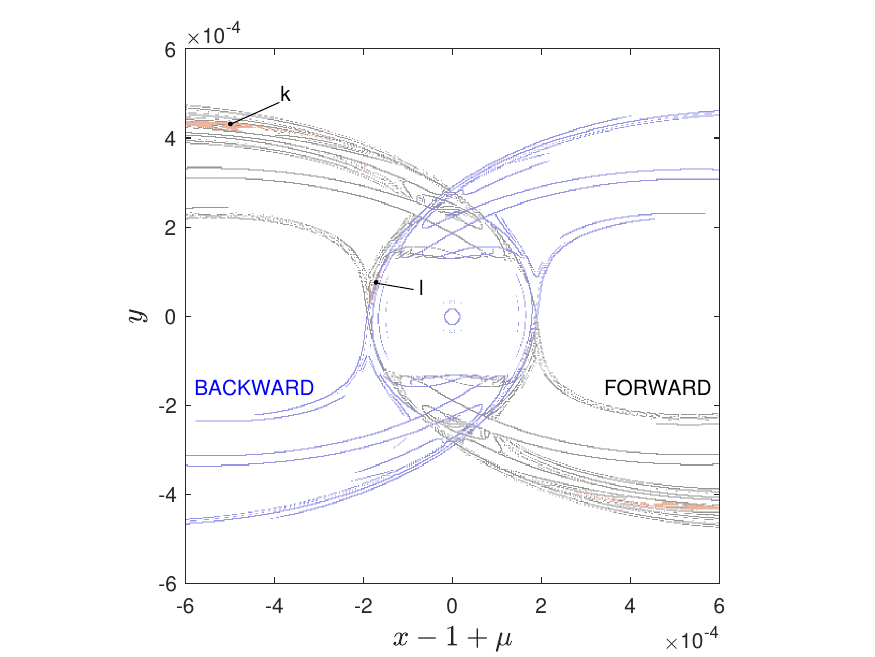}
			\caption{$\mathcal{M}(0,-\pi,3\pi)$ separatrices inspected against capture set $\mathcal{C}(-\pi,3\pi)$. Sensitivity threshold $ \sigma $ set to \num{4e-3} and \num{3e-2} for backward and forward propagations, respectively.}
			\label{fig:CaptureSets2}
	\end{subfigure}
	\caption{Inspection of separatrices extracted from $\mathcal{M}(0,f_B,f_F)$ against capture sets. Forward and backward separatrices of the \gls*{LD} field are represented as gray and blue lines, respectively.}
	\label{fig:CaptureSets}
\end{figure}

Without loss of generality, the \gls*{LD}-based approach is applied to the Sun--Mars system. The correlation between the extracted separatrices and the \gls*{WSB} is tested for several integration intervals, for both forward and backward propagations. In \Fig{fig:LD}, the \gls*{LD} scalar field computed for two distinct final anomalies is shown (see \Figsand{fig:LD1}{fig:LD3}), together with the extracted patterns overlapped to the subsets $ \mathcal{W} $, $ \mathcal{X} $, and $ \mathcal{K} $ derived for the same $f_{f}$ (see \Figsand{fig:LD2}{fig:LD4}). In both cases a good match between separatrices and boundaries of the classified regions is observed. The central green disk identifies \glspl*{IC} located inside the surface of Mars, which immediately generate crash orbits. \tr{The scalar field in \Fig{fig:LD3} is more accurate in revealing the geometrical structures that characterize the \gls*{ER3BP} dynamics, when compared to that in \Fig{fig:LD1}. The longer the finite horizon over which \glspl*{IC} are propagated, the more separatrices the field is able to reveal \cite{jimenez2009distinguished,mancho2013lagrangian,lopesino2017theoretical}. Since initial states in $\mathcal{M}(0,-\pi,0)$ are integrated over a shorter finite horizon, the field cannot reveal structures with the same level of detail as compared with $\mathcal{M}(0,0,2\pi)$, in which initial states are propagated over a longer finite horizon.} \tr{Results for the case $\mathcal{M}(0,-2\pi,0)$, here not included, reach the same accuracy of the ones for the $\mathcal{M}(0,0,2\pi)$ field. Indeed, the two \gls*{LD} fields are symmetric with respect to the $x$-axis \cite{gawlik2009lagrangian}.}

\begin{figure}[tbp]
	\centering
	\begin{subfigure}[t]{0.32\textwidth}
			\centering
			\includegraphics[width=\textwidth,clip,trim={0cm, 0.1cm, 1.2cm, 0.2cm}]{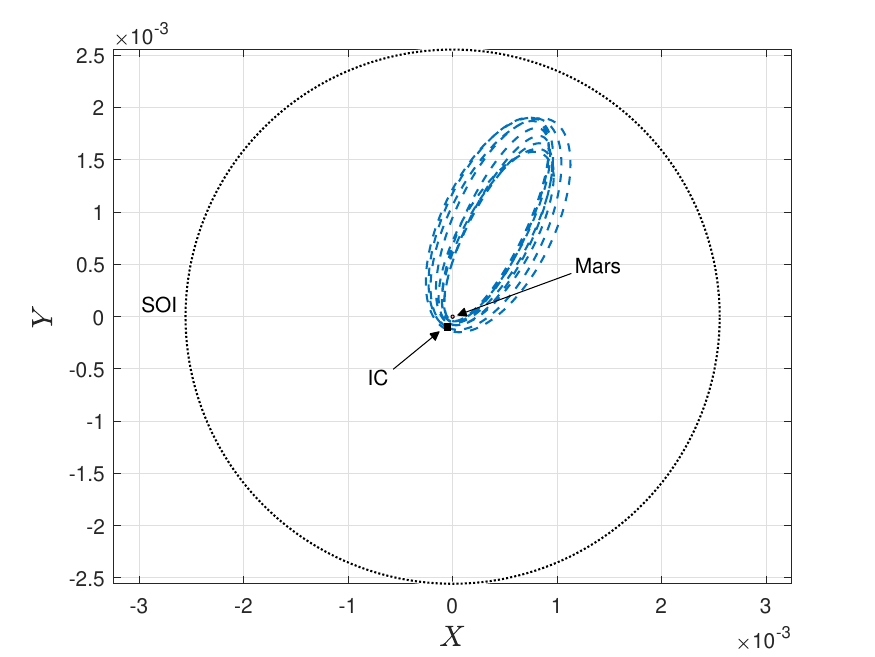}
			\caption{Orbit a.}
			\label{fig:TrajA} 
	\end{subfigure}
	\begin{subfigure}[t]{0.32\textwidth}
			\centering
			\includegraphics[width=\textwidth,clip,trim={0cm, 0.1cm, 1.2cm, 0.2cm}]{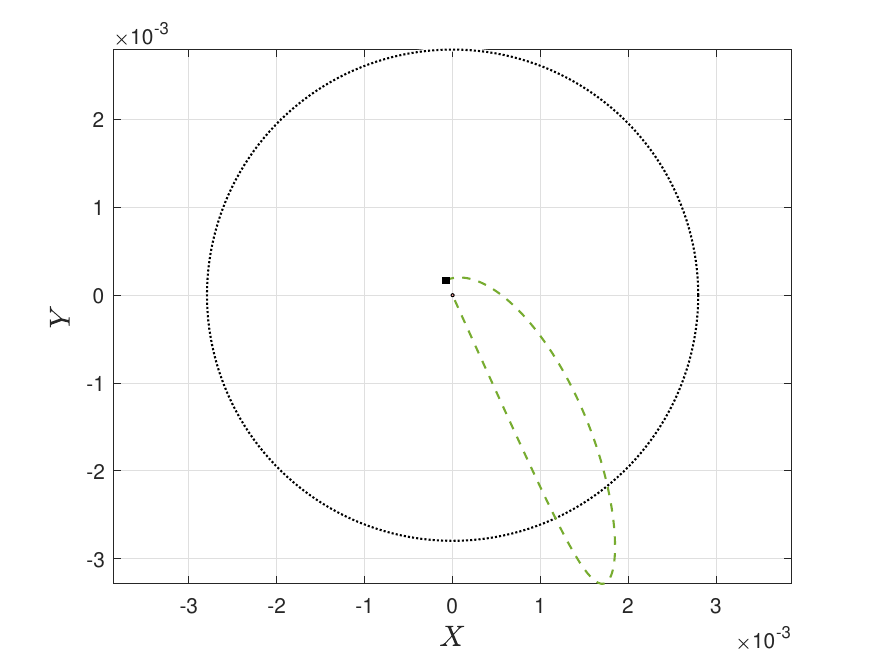}
			\caption{Orbit b.}
			\label{fig:TrajB}
	\end{subfigure}
	\begin{subfigure}[t]{0.32\textwidth}
		\centering
		\includegraphics[width=\textwidth,clip,trim={0cm, 0.1cm, 1.2cm, 0.2cm}]{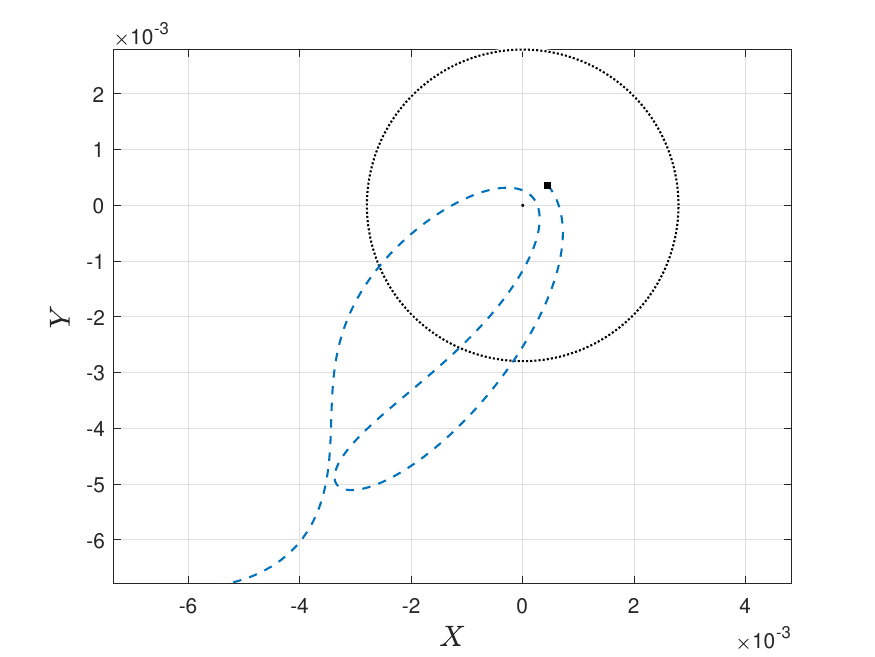}
		\caption{Orbit c.}
		\label{fig:TrajC}
	\end{subfigure}
	\begin{subfigure}[t]{0.32\textwidth}
		\centering
		\includegraphics[width=\textwidth,clip,trim={0cm, 0.1cm, 1.2cm, 0.2cm}]{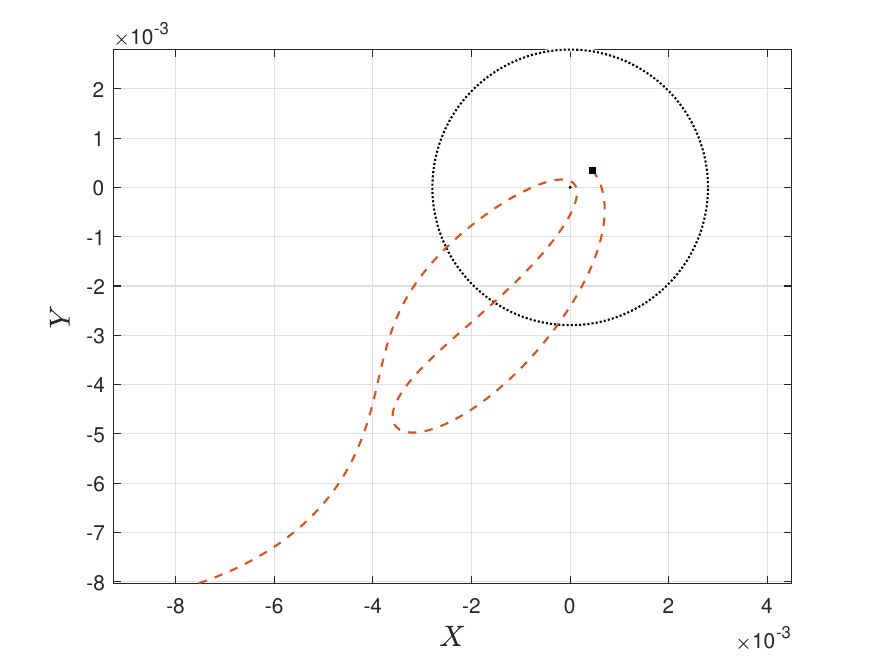}
		\caption{Orbit d.}
		\label{fig:TrajD}
	\end{subfigure}
	\begin{subfigure}[t]{0.32\textwidth}
		\centering
		\includegraphics[width=\textwidth,clip,trim={0cm, 0.1cm, 1.2cm, 0.2cm}]{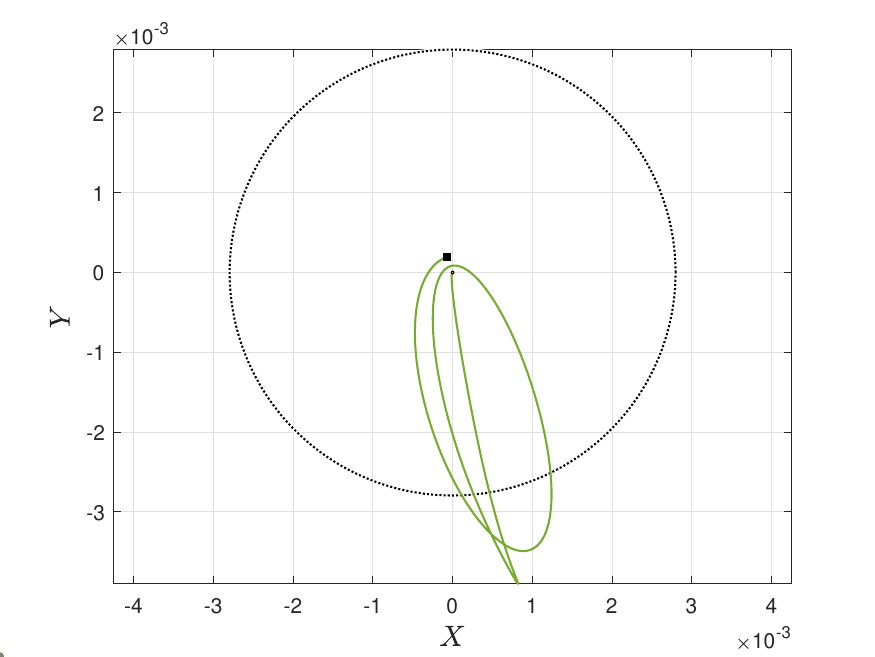}
		\caption{Orbit e.}
		\label{fig:TrajE}
	\end{subfigure}
	\begin{subfigure}[t]{0.32\textwidth}
		\centering
		\includegraphics[width=\textwidth,clip,trim={0cm, 0.1cm, 1.2cm, 0.2cm}]{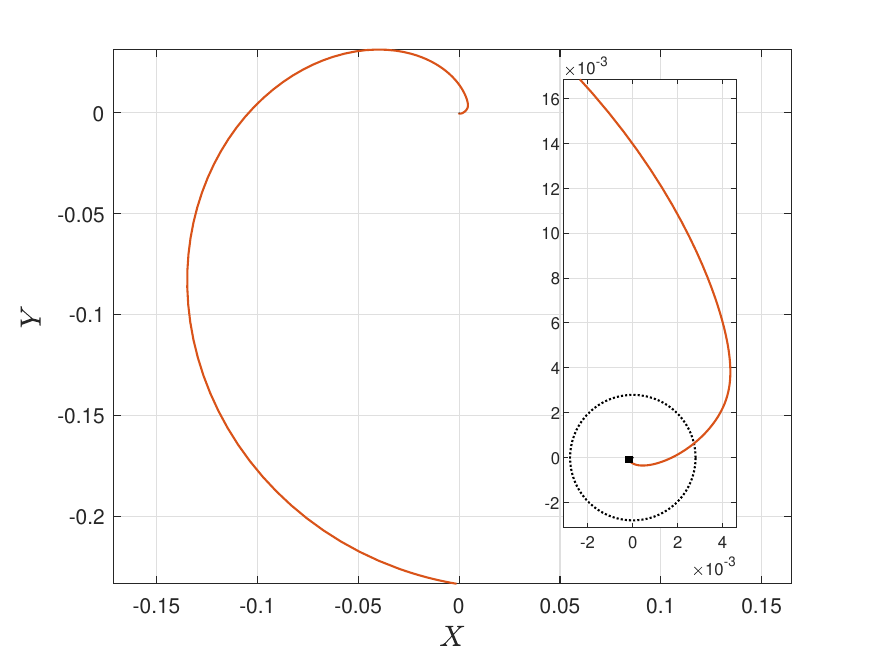}
		\caption{Orbit f.}
		\label{fig:TrajF}
	\end{subfigure}
	\begin{subfigure}[t]{0.32\textwidth}
		\centering
		\includegraphics[width=\textwidth,clip,trim={0cm, 0.1cm, 1.2cm, 0.2cm}]{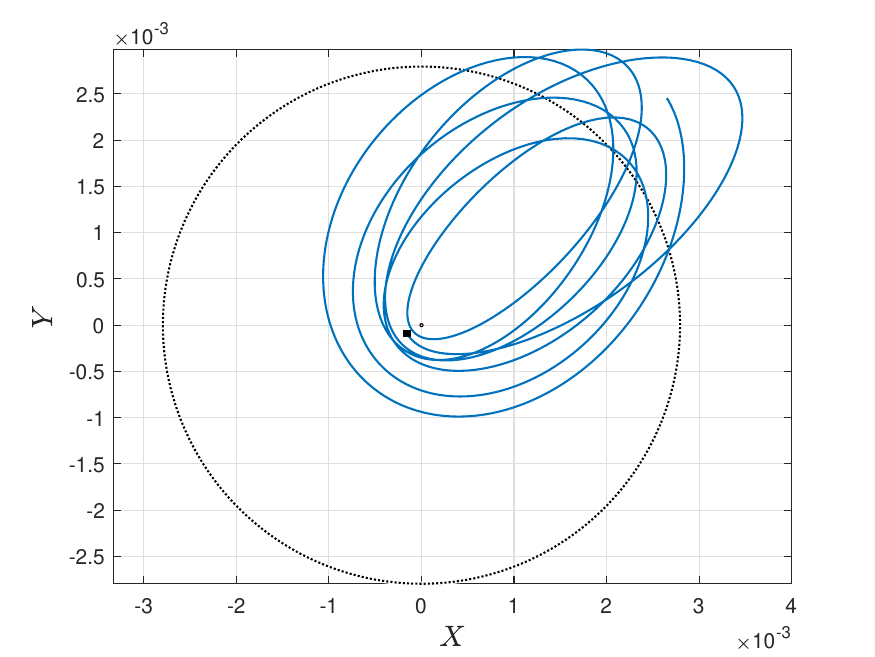}
		\caption{Orbit g.}
		\label{fig:TrajG}
	\end{subfigure}
	\begin{subfigure}[t]{0.32\textwidth}
		\centering
		\includegraphics[width=\textwidth,clip,trim={0cm, 0.1cm, 1.2cm, 0.2cm}]{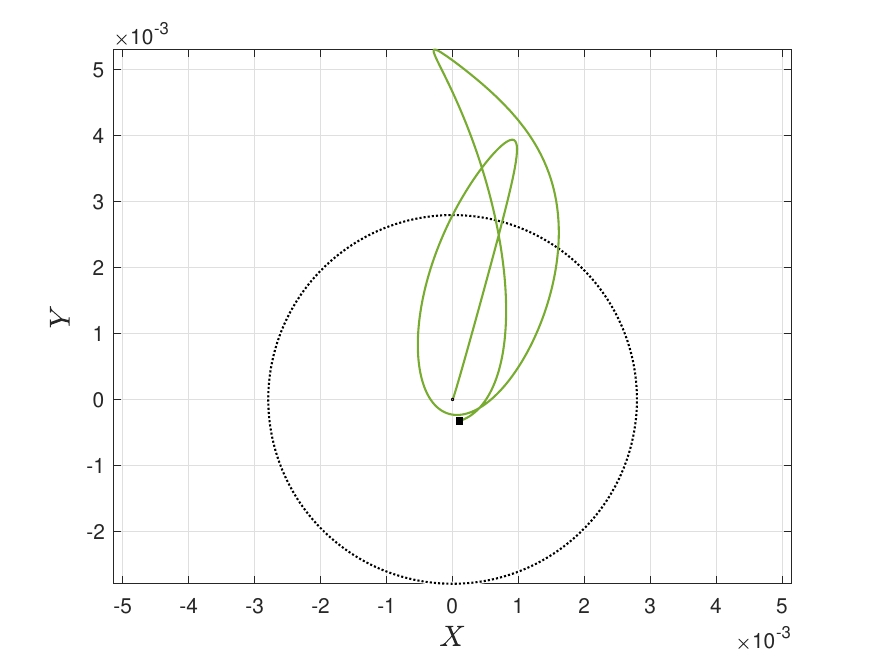}
		\caption{Orbit h.}
		\label{fig:TrajH}
	\end{subfigure}
	\begin{subfigure}[t]{0.32\textwidth}
		\centering
		\includegraphics[width=\textwidth,clip,trim={0cm, 0.1cm, 1.2cm, 0.2cm}]{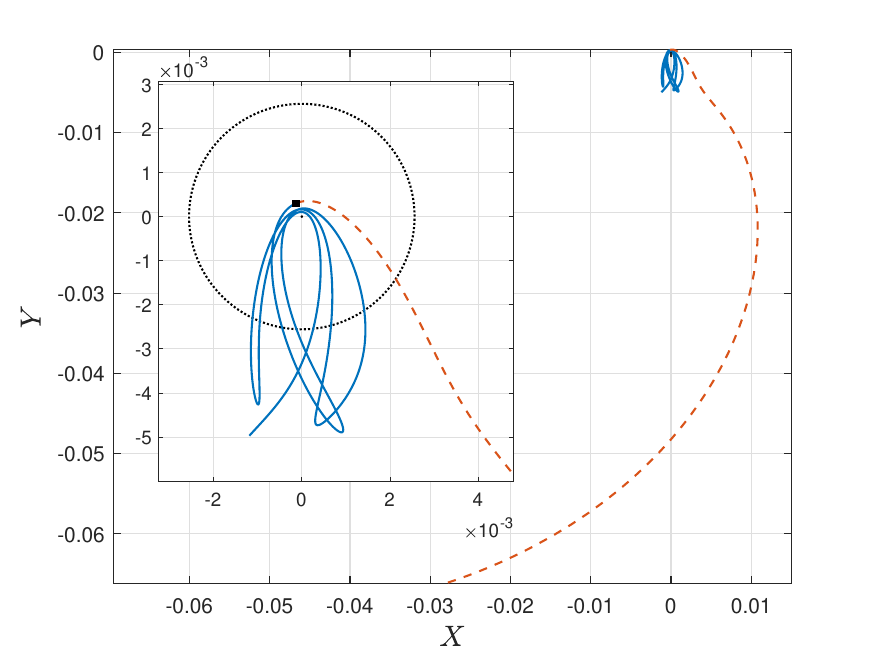}
		\caption{Orbit i.}
		\label{fig:TrajI}
	\end{subfigure}
	\begin{subfigure}[t]{0.32\textwidth}
		\centering
		\includegraphics[width=\textwidth,clip,trim={0cm, 0.1cm, 1.2cm, 0.2cm}]{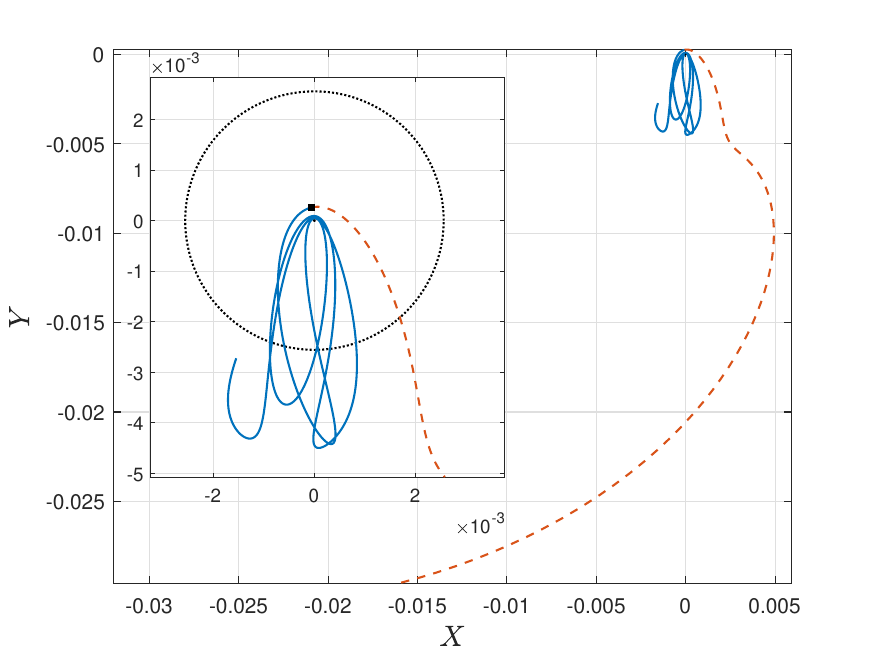}
		\caption{Orbit j.}
		\label{fig:TrajJ}
	\end{subfigure}
	\begin{subfigure}[t]{0.32\textwidth}
		\centering
		\includegraphics[width=\textwidth,clip,trim={0cm, 0.1cm, 1.2cm, 0.2cm}]{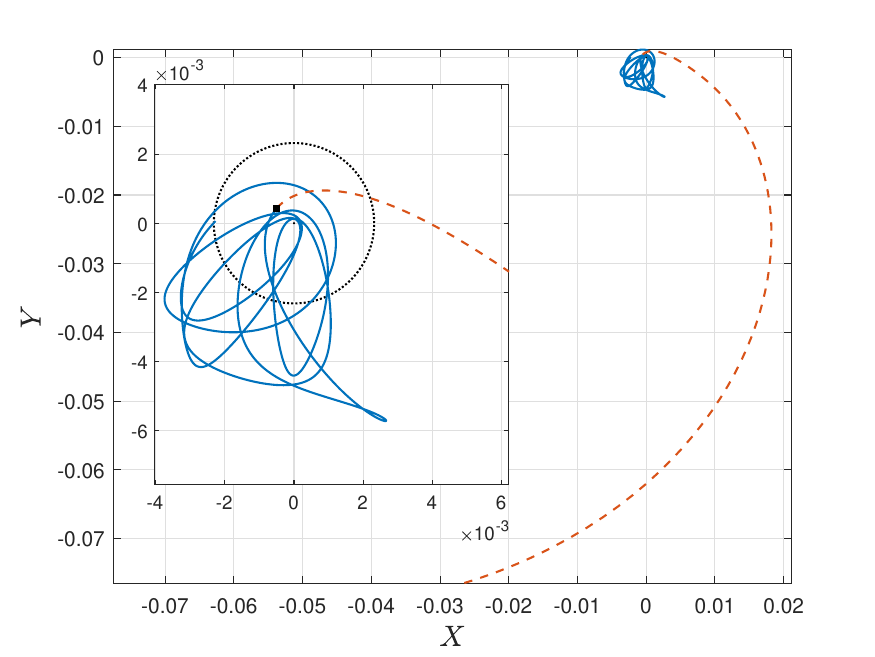}
		\caption{Orbit k.}
		\label{fig:TrajK}
	\end{subfigure}
	\begin{subfigure}[t]{0.32\textwidth}
		\centering
		\includegraphics[width=\textwidth,clip,trim={0cm, 0.1cm, 1.2cm, 0.2cm}]{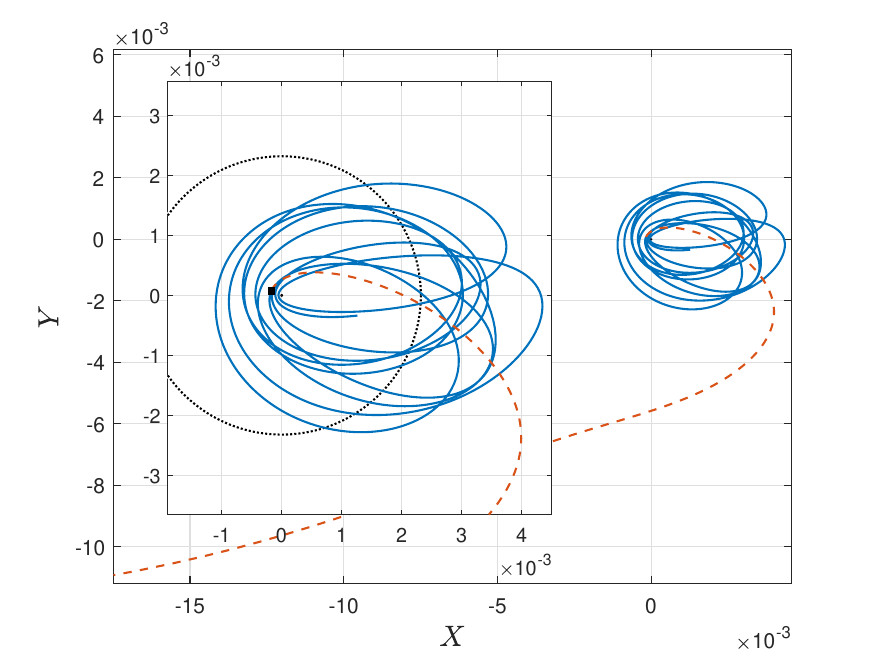}
		\caption{Orbit l.}
		\label{fig:TrajL}
	\end{subfigure}
	\caption{Sample orbits in the Mars-centered, non-rotating frame. Forward ($ f > f_{0} $) and backward ($ f < f_{0} $) legs plotted as solid and dashed lines, respectively; Mars' \gls*{SOI} represented as a black dotted circumference. \glspl*{IC} indicated with square markers. Coloring identifies the subset the orbits belong to, according to the color code used in \Fig{fig:LD}. \glspl*{IC} collected in \Tab{tab:ics}.}
	\label{fig:Trajectories}
\end{figure}

For small values of $ f_{f} $, the matching presents some inconsistencies that are intrinsic to the \gls*{LD} definition. In fact, \gls*{LD} reveals patterns if \glspl*{IC} are integrated long enough for dynamical divergences between orbits to be manifested \cite{mancho2013lagrangian}. Consequently, the classification of the phase space according to the definition of particle stability provided in \Sec{sec:stability-def} may be inconsistent with some regions featured by the \gls*{LD} scalar field if the trajectories are not sufficiently divergent to feature singular structures in the field \cite{mancho2013lagrangian}. The latter is particularly true for short integration intervals as observed in \Fig{fig:LD2}. For instance, \glspl*{IC} `c' and `d' in \Fig{fig:LD2} are classified into two different sets, still their dynamical behavior is very similar as shown by their orbits in \Figsand{fig:TrajC}{fig:TrajD}. Indeed, for a slightly larger integration interval both orbits escape from Mars.

Remarkably, \glspl*{LD} detect divergence (forward propagation) and attraction (backward propagation) in the dynamical behavior even in areas classified in the same way according to our particle stability definition. To illustrate this concept, two grid points can evolve both in crash orbits, nonetheless their trajectories could be strongly different, as well as their impact epochs. For example, samples `e' and `h' in \Fig{fig:LD4} belong to two distinct regions of the same crash set $ \mathcal{K}(2\pi) $, therefore they both impact with Mars. However, they exhibit dissimilar trajectories (see \Figsand{fig:TrajE}{fig:TrajH}). They impact from different directions, and orbit `h' reverses its angular momentum with respect to Mars much earlier than orbit `e'. 

Patterns ruling particles transport in both true anomaly directions are revealed combining the \gls*{LD} structures propagated forwards and backwards \cite{mancho2013lagrangian}. The correlation of the two capture sets $ \mathcal{C}(-\pi,3\pi/2) $ and $\mathcal{C}(-\pi,3\pi)$ with the separatrices extracted from $\mathcal{M}(0,-\pi,3\pi/2)$ and $\mathcal{M}(0,-\pi,3\pi)$ fields, respectively, is presented in \Fig{fig:CaptureSets}. Results show that some of the areas in the phase space enclosed by \gls*{LD} separatrices appear to be capture sets. Based on the outcome of the validation procedure, the devised methodology of computing the \gls*{LD} field and extracting the dynamics separatrices has been proven successful. 

Referring to \Fig{fig:CaptureSets}, the \gls*{LD} approach omits the dynamical behavior featured by the highlighted numerous regions, therefore a classification technique is still required to discern which areas are actually capture sets. A viable strategy to overcome the aforementioned limitation is proposed for practical design of \gls*{BC} orbits. By sampling an individual \gls*{IC} for each identified region and classifying its orbit, all areas in the phase space can be easily categorized either as $ \mathcal{W} $, $ \mathcal{X} $, or $ \mathcal{K} $ subsets according to the particle stability definition given in \Sec{sec:stability-def}.

The exact \glspl*{IC} sampled from \Figsand{fig:LD}{fig:CaptureSets} are collected in \Tab{tab:ics}. Their orbits expressed in the Mars-centered, non-rotating frame, oriented as the synodic frame at $f_{0}$, and with Cartesian coordinates $X$ and $Y$ are plotted in \Fig{fig:Trajectories}. Compared to similar \gls*{BC} orbits found in the literature \cite{hyeraci2010method,luo2014constructing}, the weakly stable trajectories shown in \Fig{fig:Trajectories} do not fully complete the last revolution about Mars due to the dropping of the usual geometrical constraint on the revolutions number. Nevertheless, they grant temporary capture at least over the finite horizon specified by the integration interval.

\begin{table}[pos=tbp,width=0.9\textwidth,align=\centering]
	\caption{Initial conditions of sample orbits.}
	\label{tab:ics}
	\begin{tabular*}{0.9\textwidth}{@{}CRRRRL@{}}
		\toprule
		\multirow{2}{*}{\textbf{Orbit}} & \multicolumn{4}{c}{\textbf{Initial condition at} $ f_{0} = 0 $} & \multirow{2}{*}{\textbf{Set}} \\
		\cmidrule(lr){2-5}
		& \multicolumn{1}{c}{$ X_{0} = x_{0}-1+\mu $} & \multicolumn{1}{c}{$ Y_{0} = y_{0} $} & \multicolumn{1}{c}{$ {{x}'}_{0} $} & \multicolumn{1}{c}{$ {{y}'}_{0} $} & \\
		\midrule
		a & \num{-5.170000e-5} & \num{-1.000000e-4} & \num{6.258637e-02} &   \num{-3.235715e-02} & $ \mathcal{W}(-\pi) $ \\
		b & \num{-7.575000e-05} & \num{1.695000e-04} & \num{-4.999940e-02} & \num{-2.234486e-02} & $ \mathcal{K}(-\pi) $ \\
		c & \num{4.509000e-04} & \num{3.621000e-04} & \num{-1.913330e-02} & \num{2.382548e-02} & $ \mathcal{W}(-\pi) $ \\
		d & \num{4.533000e-04} & \num{3.475000e-04} & \num{-1.871302e-02} & \num{2.441039e-02} & $ \mathcal{X}(-\pi) $ \\
		e & \num{-7.094000e-05} & \num{1.960000e-04} & \num{-4.856863e-02} & \num{-1.757887e-02} & $ \mathcal{K}(2\pi) $ \\
		f & \num{-1.719000e-04} & \num{-9.739000e-05} & \num{2.616042e-02} & \num{-4.617492e-02} & $ \mathcal{X}(2\pi) $ \\
		g & \num{-1.551000e-04} & \num{-9.239000e-05} & \num{2.842583e-02} & \num{-4.771994e-02} & $ \mathcal{W}(2\pi) $ \\
		h & \num{1.094000e-04} & \num{-3.258000e-04} & \num{3.796152e-02} & \num{1.274705e-02} & $ \mathcal{K}(2\pi) $ \\
		i & \num{-1.286000e-04} & \num{3.018000e-04} & \num{-3.772815e-02} & \num{-1.607634e-02} & $ \mathcal{C}(-\pi,3\pi/2) $ \\
		j & \num{-6.373000e-05} & \num{2.585000e-04} & \num{-4.429485e-02} & \num{-1.092035e-02} & $ \mathcal{C}(-\pi,3\pi/2) $ \\
		k & \num{-4.990000e-04} & \num{4.317000e-04} & \num{-1.863920e-02} & \num{-2.154496e-02} & $ \mathcal{C}(-\pi,3\pi) $ \\
		l & \num{-1.719000e-04} & \num{7.575000e-05} & \num{-2.195327e-02} & \num{-4.981872e-02} & $ \mathcal{C}(-\pi,3\pi) $ \\
		\bottomrule
	\end{tabular*}
\end{table}

\section{Conclusion} \label{sec:conclusion}
This study looks into the effectiveness of \acrlongpl*{LD} in revealing phase space organizing structures when facing the \acrlong*{BC} phenomenon, so investigating their capability in highlighting the \acrlong*{WSB}. They successfully reveal the geometrical structures governing the transport mechanisms that are then extracted through an edge detection algorithm based on the Roberts' operator to approximate the gradient of the \acrlong*{LD} scalar fields. The extracted separatrices effectively distinguish regions with different dynamical behavior. The detected patterns are in good agreement with the \acrlong*{WSB} computed on the same integration interval. \acrlongpl*{LD} proved to be an intuitive, easy to implement, and convenient tool for designing \acrlong*{BC} orbits. Without any a priory knowledge, \acrlong*{LD} patterns yield a consistent match with the \acrlong*{WSB} and the associated stable sets. The proposed methodology supports the design of \acrlong*{BC} orbits, enriching the dynamics knowledge in the proximity of the target planet. Furthermore, the technique can be successfully applied to arbitrary non-autonomous, more representative, astrodynamics models without any restrictions.

\printcredits

\section*{Declaration of competing interest}

The authors declare that they have no known competing financial interests or personal relationships that could have appeared to influence the work reported in this paper.

\section*{Acknowledgment}

G.M. and F.T. would like to acknowledge the \gls*{ERC} since part of this work has received funding from the \gls*{ERC} under the European Union’s Horizon 2020 research and innovation programme (Grant Agreement No.\,864697).

\bibliographystyle{cas-model2-names}

\bibliography{references}

\begin{thebibliography}{37}
\expandafter\ifx\csname natexlab\endcsname\relax\def\natexlab#1{#1}\fi
\providecommand{\url}[1]{\texttt{#1}}
\providecommand{\href}[2]{#2}
\providecommand{\path}[1]{#1}
\providecommand{\DOIprefix}{doi:}
\providecommand{\ArXivprefix}{arXiv:}
\providecommand{\URLprefix}{URL: }
\providecommand{\Pubmedprefix}{pmid:}
\providecommand{\doi}[1]{\href{http://dx.doi.org/#1}{\path{#1}}}
\providecommand{\Pubmed}[1]{\href{pmid:#1}{\path{#1}}}
\providecommand{\bibinfo}[2]{#2}
\ifx\xfnm\relax \def\xfnm[#1]{\unskip,\space#1}\fi
\bibitem[{Belbruno(1987)}]{belbruno1987lunar}
\bibinfo{author}{Belbruno, E.}, \bibinfo{year}{1987}.
\newblock \bibinfo{title}{Lunar capture orbits, a method of constructing
  {Earth} {Moon} trajectories and the lunar {GAS} mission}, in:
  \bibinfo{booktitle}{19th International Electric Propulsion Conference},
  \bibinfo{publisher}{American Institute of Aeronautics and Astronautics}. pp.
  \bibinfo{pages}{1--9}.
\newblock \DOIprefix\doi{10.2514/6.1987-1054}.
\bibitem[{Belbruno(2004)}]{belbruno2004capture}
\bibinfo{author}{Belbruno, E.}, \bibinfo{year}{2004}.
\newblock \bibinfo{title}{Capture Dynamics and Chaotic Motions in Celestial
  Mechanics}.
\newblock \bibinfo{publisher}{Princeton University Press}.
\newblock \DOIprefix\doi{10.1515/9780691186436}.
\bibitem[{Belbruno and Carrico(2000)}]{belbruno2000calculation}
\bibinfo{author}{Belbruno, E.}, \bibinfo{author}{Carrico, J.},
  \bibinfo{year}{2000}.
\newblock \bibinfo{title}{Calculation of weak stability boundary ballistic
  lunar transfer trajectories}, in: \bibinfo{booktitle}{Astrodynamics
  Specialist Conference}, pp. \bibinfo{pages}{4142, 262--271}.
\newblock \DOIprefix\doi{10.2514/6.2000-4142}.
\bibitem[{Belbruno et~al.(2010)Belbruno, Gidea and Topputo}]{belbruno2010weak}
\bibinfo{author}{Belbruno, E.}, \bibinfo{author}{Gidea, M.},
  \bibinfo{author}{Topputo, F.}, \bibinfo{year}{2010}.
\newblock \bibinfo{title}{Weak stability boundary and invariant manifolds}.
\newblock \bibinfo{journal}{SIAM Journal on Applied Dynamical Systems}
  \bibinfo{volume}{9}, \bibinfo{pages}{1061--1089}.
\newblock \DOIprefix\doi{10.1137/090780638}.
\bibitem[{Belbruno and Miller(1990)}]{belbruno1990ballistic}
\bibinfo{author}{Belbruno, E.}, \bibinfo{author}{Miller, J.},
  \bibinfo{year}{1990}.
\newblock \bibinfo{title}{A ballistic lunar capture trajectory for the
  {Japanese} spacecraft {Hiten}}.
\newblock \bibinfo{type}{Technical Report}. Jet Propulsion Laboratory.
\newblock \bibinfo{note}{{IOM} 312/90.4-1731-EAB}.
\bibitem[{Belbruno and Miller(1993)}]{belbruno1993sun}
\bibinfo{author}{Belbruno, E.}, \bibinfo{author}{Miller, J.},
  \bibinfo{year}{1993}.
\newblock \bibinfo{title}{{Sun}-perturbed {Earth}-to-{Moon} transfers with
  ballistic capture}.
\newblock \bibinfo{journal}{Journal of Guidance, Control, and Dynamics}
  \bibinfo{volume}{16}, \bibinfo{pages}{770--775}.
\newblock \DOIprefix\doi{10.2514/3.21079}.
\bibitem[{Bernardini et~al.(2022)Bernardini, Merisio, Losacco, Raffa, Armellin,
  Baresi, Lizy-Destrez, Topputo, Canalias et~al.}]{bernardini2022exploiting}
\bibinfo{author}{Bernardini, N.}, \bibinfo{author}{Merisio, G.},
  \bibinfo{author}{Losacco, M.}, \bibinfo{author}{Raffa, S.},
  \bibinfo{author}{Armellin, R.}, \bibinfo{author}{Baresi, N.},
  \bibinfo{author}{Lizy-Destrez, S.}, \bibinfo{author}{Topputo, F.},
  \bibinfo{author}{Canalias, E.}, et~al., \bibinfo{year}{2022}.
\newblock \bibinfo{title}{\tr{Exploiting coherent patterns for the analysis of
  qualitative motion and the design of bounded orbits around small bodies}},
  in: \bibinfo{booktitle}{73rd International Astronautical Congress (IAC
  2022)}, pp. \bibinfo{pages}{1--16}.
\bibitem[{Caleb et~al.(2022)Caleb, Merisio, Di~Lizia and
  Topputo}]{caleb2022stable}
\bibinfo{author}{Caleb, T.}, \bibinfo{author}{Merisio, G.},
  \bibinfo{author}{Di~Lizia, P.}, \bibinfo{author}{Topputo, F.},
  \bibinfo{year}{2022}.
\newblock \bibinfo{title}{\tr{Stable sets mapping with {Taylor} differential
  algebra with application to ballistic capture orbits around {Mars}}}.
\newblock \bibinfo{journal}{Celestial Mechanics and Dynamical Astronomy}
  \bibinfo{volume}{134}, \bibinfo{pages}{39}.
\newblock \DOIprefix\doi{10.1007/s10569-022-10090-8}.
\bibitem[{Campagnola and Russell(2010)}]{campagnola2010endgame}
\bibinfo{author}{Campagnola, S.}, \bibinfo{author}{Russell, R.P.},
  \bibinfo{year}{2010}.
\newblock \bibinfo{title}{\tr{Endgame problem part 2: Multibody technique and
  the {Tisserand}-{Poincare} graph}}.
\newblock \bibinfo{journal}{Journal of Guidance, Control, and Dynamics}
  \bibinfo{volume}{33}, \bibinfo{pages}{476--486}.
\newblock \DOIprefix\doi{10.2514/1.44290}.
\bibitem[{Carletta et~al.(2019)Carletta, Pontani and
  Teofilatto}]{carletta2019design}
\bibinfo{author}{Carletta, S.}, \bibinfo{author}{Pontani, M.},
  \bibinfo{author}{Teofilatto, P.}, \bibinfo{year}{2019}.
\newblock \bibinfo{title}{\tr{Design of low-energy capture trajectories in the
  elliptic restricted four-body problem}}, in: \bibinfo{booktitle}{70th
  International Astronautical Congress (IAC 2019)}, pp.
  \bibinfo{pages}{21--25}.
\bibitem[{Circi and Teofilatto(2001)}]{circi2001dynamics}
\bibinfo{author}{Circi, C.}, \bibinfo{author}{Teofilatto, P.},
  \bibinfo{year}{2001}.
\newblock \bibinfo{title}{On the dynamics of weak stability boundary lunar
  transfers}.
\newblock \bibinfo{journal}{Celestial Mechanics and Dynamical Astronomy}
  \bibinfo{volume}{79}, \bibinfo{pages}{41--72}.
\newblock \DOIprefix\doi{10.1023/A:1011153610564}.
\bibitem[{Conley(1968)}]{conley1968low}
\bibinfo{author}{Conley, C.C.}, \bibinfo{year}{1968}.
\newblock \bibinfo{title}{\tr{Low energy transit orbits in the restricted
  three-body problems}}.
\newblock \bibinfo{journal}{SIAM Journal on Applied Mathematics}
  \bibinfo{volume}{16}, \bibinfo{pages}{732--746}.
\newblock \DOIprefix\doi{10.1137/0116060}.
\bibitem[{Conley(1969)}]{conley1969ultimate}
\bibinfo{author}{Conley, C.C.}, \bibinfo{year}{1969}.
\newblock \bibinfo{title}{\tr{On the ultimate behavior of orbits with respect
  to an unstable critical point I. Oscillating, asymptotic, and capture
  orbits}}.
\newblock \bibinfo{journal}{Journal of Differential Equations}
  \bibinfo{volume}{5}, \bibinfo{pages}{136--158}.
\newblock \DOIprefix\doi{10.1016/0022-0396(69)90108-9}.
\bibitem[{Davis(1975)}]{davis1975survey}
\bibinfo{author}{Davis, L.S.}, \bibinfo{year}{1975}.
\newblock \bibinfo{title}{A survey of edge detection techniques}.
\newblock \bibinfo{journal}{Computer Graphics and Image Processing}
  \bibinfo{volume}{4}, \bibinfo{pages}{248--270}.
\newblock \DOIprefix\doi{10.1016/0146-664X(75)90012-X}.
\bibitem[{Garc{\'i}a and G{\'o}mez(2007)}]{garcia2007note}
\bibinfo{author}{Garc{\'i}a, F.}, \bibinfo{author}{G{\'o}mez, G.},
  \bibinfo{year}{2007}.
\newblock \bibinfo{title}{A note on weak stability boundaries}.
\newblock \bibinfo{journal}{Celestial Mechanics and Dynamical Astronomy}
  \bibinfo{volume}{97}, \bibinfo{pages}{87--100}.
\newblock \DOIprefix\doi{10.1007/s10569-006-9053-6}.
\bibitem[{Gawlik et~al.(2009)Gawlik, Marsden, Du~Toit and
  Campagnola}]{gawlik2009lagrangian}
\bibinfo{author}{Gawlik, E.S.}, \bibinfo{author}{Marsden, J.E.},
  \bibinfo{author}{Du~Toit, P.C.}, \bibinfo{author}{Campagnola, S.},
  \bibinfo{year}{2009}.
\newblock \bibinfo{title}{\tr{{Lagrangian} coherent structures in the planar
  elliptic restricted three-body problem}}.
\newblock \bibinfo{journal}{Celestial Mechanics and Dynamical Astronomy}
  \bibinfo{volume}{103}, \bibinfo{pages}{227--249}.
\bibitem[{Haller(2011)}]{haller2011variational}
\bibinfo{author}{Haller, G.}, \bibinfo{year}{2011}.
\newblock \bibinfo{title}{A variational theory of hyperbolic {Lagrangian}
  coherent structures}.
\newblock \bibinfo{journal}{Physica D: Nonlinear Phenomena}
  \bibinfo{volume}{240}, \bibinfo{pages}{574--598}.
\newblock \DOIprefix\doi{10.1016/j.physd.2010.11.010}.
\bibitem[{Haller(2015)}]{haller2015lagrangian}
\bibinfo{author}{Haller, G.}, \bibinfo{year}{2015}.
\newblock \bibinfo{title}{{Lagrangian} coherent structures}.
\newblock \bibinfo{journal}{Annual Review of Fluid Mechanics}
  \bibinfo{volume}{47}, \bibinfo{pages}{137--162}.
\newblock \DOIprefix\doi{10.1146/annurev-fluid-010313-141322}.
\bibitem[{Hyeraci and Topputo(2010)}]{hyeraci2010method}
\bibinfo{author}{Hyeraci, N.}, \bibinfo{author}{Topputo, F.},
  \bibinfo{year}{2010}.
\newblock \bibinfo{title}{Method to design ballistic capture in the elliptic
  restricted three-body problem}.
\newblock \bibinfo{journal}{Journal of Guidance, Control, and Dynamics}
  \bibinfo{volume}{33}, \bibinfo{pages}{1814--1823}.
\newblock \DOIprefix\doi{10.2514/1.49263}.
\bibitem[{Jim{\'e}nez~Madrid and Mancho(2009)}]{jimenez2009distinguished}
\bibinfo{author}{Jim{\'e}nez~Madrid, J.A.}, \bibinfo{author}{Mancho, A.M.},
  \bibinfo{year}{2009}.
\newblock \bibinfo{title}{Distinguished trajectories in time dependent vector
  fields}.
\newblock \bibinfo{journal}{Chaos: An Interdisciplinary Journal of Nonlinear
  Science} \bibinfo{volume}{19}, \bibinfo{pages}{013111, 1--19}.
\newblock \DOIprefix\doi{10.1063/1.3056050}.
\bibitem[{Lopesino et~al.(2017)Lopesino, Balibrea-Iniesta, Garc{\'\i}a-Garrido,
  Wiggins and Mancho}]{lopesino2017theoretical}
\bibinfo{author}{Lopesino, C.}, \bibinfo{author}{Balibrea-Iniesta, F.},
  \bibinfo{author}{Garc{\'\i}a-Garrido, V.J.}, \bibinfo{author}{Wiggins, S.},
  \bibinfo{author}{Mancho, A.M.}, \bibinfo{year}{2017}.
\newblock \bibinfo{title}{A theoretical framework for {Lagrangian}
  descriptors}.
\newblock \bibinfo{journal}{International Journal of Bifurcation and Chaos}
  \bibinfo{volume}{27}, \bibinfo{pages}{1730001, 1--25}.
\newblock \DOIprefix\doi{10.1142/s0218127417300014}.
\bibitem[{{Luo, Z.-F.} and Topputo(2015)}]{luo2015analysis}
\bibinfo{author}{{Luo, Z.-F.}}, \bibinfo{author}{Topputo, F.},
  \bibinfo{year}{2015}.
\newblock \bibinfo{title}{Analysis of ballistic capture in {Sun}--planet
  models}.
\newblock \bibinfo{journal}{Advances in Space Research} \bibinfo{volume}{56},
  \bibinfo{pages}{1030--1041}.
\newblock \DOIprefix\doi{10.1016/j.asr.2015.05.042}.
\bibitem[{{Luo, Z.-F.} et~al.(2014){Luo, Z.-F.}, Topputo, Bernelli-Zazzera and
  {Tang, G.-J.}}]{luo2014constructing}
\bibinfo{author}{{Luo, Z.-F.}}, \bibinfo{author}{Topputo, F.},
  \bibinfo{author}{Bernelli-Zazzera, F.}, \bibinfo{author}{{Tang, G.-J.}},
  \bibinfo{year}{2014}.
\newblock \bibinfo{title}{Constructing ballistic capture orbits in the real
  solar system model}.
\newblock \bibinfo{journal}{Celestial Mechanics and Dynamical Astronomy}
  \bibinfo{volume}{120}, \bibinfo{pages}{433--450}.
\newblock \DOIprefix\doi{10.1007/s10569-014-9580-5}.
\bibitem[{Mancho et~al.(2013)Mancho, Wiggins, Curbelo and
  Mendoza}]{mancho2013lagrangian}
\bibinfo{author}{Mancho, A.M.}, \bibinfo{author}{Wiggins, S.},
  \bibinfo{author}{Curbelo, J.}, \bibinfo{author}{Mendoza, C.},
  \bibinfo{year}{2013}.
\newblock \bibinfo{title}{{Lagrangian} descriptors: A method for revealing
  phase space structures of general time dependent dynamical systems}.
\newblock \bibinfo{journal}{Communications in Nonlinear Science and Numerical
  Simulation} \bibinfo{volume}{18}, \bibinfo{pages}{3530--3557}.
\newblock \DOIprefix\doi{10.1016/j.cnsns.2013.05.002}.
\bibitem[{Manzi and Topputo(2021)}]{manzi2021flow}
\bibinfo{author}{Manzi, M.}, \bibinfo{author}{Topputo, F.},
  \bibinfo{year}{2021}.
\newblock \bibinfo{title}{A flow-informed strategy for ballistic capture orbit
  generation}.
\newblock \bibinfo{journal}{Celestial Mechanics and Dynamical Astronomy}
  \bibinfo{volume}{133}, \bibinfo{pages}{1--16}.
\newblock \DOIprefix\doi{10.1007/s10569-021-10048-2}.
\bibitem[{Merisio(2023)}]{merisio2023engineering}
\bibinfo{author}{Merisio, G.}, \bibinfo{year}{2023}.
\newblock \bibinfo{title}{\tr{Engineering ballistic capture for autonomous
  interplanetary spacecraft with limited onboard resources ({PhD} Thesis)}}.
\newblock \bibinfo{publisher}{Politecnico di Milano, Milan, Italy}.
\newblock \URLprefix \url{http://hdl.handle.net/10589/196152}.
\bibitem[{Montenbruck and Gill(2000)}]{montenbruck2000satellite}
\bibinfo{author}{Montenbruck, O.}, \bibinfo{author}{Gill, E.},
  \bibinfo{year}{2000}.
\newblock \bibinfo{title}{Satellite Orbits Models, Methods and Applications}.
\newblock \bibinfo{publisher}{Springer}.
\newblock \DOIprefix\doi{10.1007/978-3-642-58351-3}.
\bibitem[{Raffa et~al.(2023)Raffa, Merisio and Topputo}]{raffa2023finding}
\bibinfo{author}{Raffa, S.}, \bibinfo{author}{Merisio, G.},
  \bibinfo{author}{Topputo, F.}, \bibinfo{year}{2023}.
\newblock \bibinfo{title}{\tr{Finding regions of bounded motion in binary
  asteroid environment using {Lagrangian} descriptors}}.
\newblock \bibinfo{journal}{Communications in Nonlinear Science and Numerical
  Simulation} ,
  \bibinfo{pages}{107198}\DOIprefix\doi{10.1016/j.cnsns.2023.107198}.
\bibitem[{Restrepo and Russell(2017)}]{restrepo2017patched}
\bibinfo{author}{Restrepo, R.L.}, \bibinfo{author}{Russell, R.P.},
  \bibinfo{year}{2017}.
\newblock \bibinfo{title}{\tr{Patched periodic orbits: A systematic strategy
  for low energy transfer design}}, in: \bibinfo{booktitle}{AAS/AIAA
  Astrodynamics Specialist Conference}, pp. \bibinfo{pages}{17--695}.
\bibitem[{Roberts(1963)}]{roberts1963machine}
\bibinfo{author}{Roberts, L.G.}, \bibinfo{year}{1963}.
\newblock \bibinfo{title}{Machine perception of three-dimensional solids ({PhD}
  Thesis)}.
\newblock \bibinfo{publisher}{Massachusetts Institute of Technology}.
\bibitem[{Sousa~Silva and Terra(2012)}]{silva2012applicability}
\bibinfo{author}{Sousa~Silva, P.}, \bibinfo{author}{Terra, M.},
  \bibinfo{year}{2012}.
\newblock \bibinfo{title}{Applicability and dynamical characterization of the
  associated sets of the algorithmic weak stability boundary in the lunar
  sphere of influence}.
\newblock \bibinfo{journal}{Celestial Mechanics and Dynamical Astronomy}
  \bibinfo{volume}{113}, \bibinfo{pages}{141--168}.
\newblock \DOIprefix\doi{10.1007/s10569-012-9409-z}.
\bibitem[{Topputo and Belbruno(2009)}]{topputo2009computation}
\bibinfo{author}{Topputo, F.}, \bibinfo{author}{Belbruno, E.},
  \bibinfo{year}{2009}.
\newblock \bibinfo{title}{Computation of weak stability boundaries:
  {Sun}--{Jupiter} system}.
\newblock \bibinfo{journal}{Celestial Mechanics and Dynamical Astronomy}
  \bibinfo{volume}{105}, \bibinfo{pages}{3--17}.
\newblock \DOIprefix\doi{10.1007/s10569-009-9222-5}.
\bibitem[{Topputo and Belbruno(2015)}]{topputo2015earth}
\bibinfo{author}{Topputo, F.}, \bibinfo{author}{Belbruno, E.},
  \bibinfo{year}{2015}.
\newblock \bibinfo{title}{{Earth}--{Mars} transfers with ballistic capture}.
\newblock \bibinfo{journal}{Celestial Mechanics and Dynamical Astronomy}
  \bibinfo{volume}{121}, \bibinfo{pages}{329--346}.
\newblock \DOIprefix\doi{10.1007/s10569-015-9605-8}.
\bibitem[{Topputo et~al.(2005)Topputo, Vasile and
  Bernelli-Zazzera}]{topputo2005low}
\bibinfo{author}{Topputo, F.}, \bibinfo{author}{Vasile, M.},
  \bibinfo{author}{Bernelli-Zazzera, F.}, \bibinfo{year}{2005}.
\newblock \bibinfo{title}{Low energy interplanetary transfers exploiting
  invariant manifolds of the restricted three-body problem}.
\newblock \bibinfo{journal}{The Journal of the Astronautical Sciences}
  \bibinfo{volume}{53}, \bibinfo{pages}{353--372}.
\newblock \DOIprefix\doi{10.1007/BF03546358}.
\bibitem[{Tyler and Wittig(2022)}]{tyler2022improved}
\bibinfo{author}{Tyler, J.}, \bibinfo{author}{Wittig, A.},
  \bibinfo{year}{2022}.
\newblock \bibinfo{title}{\tr{An improved numerical method for hyperbolic
  {Lagrangian} coherent structures using differential algebra}}.
\newblock \bibinfo{journal}{Journal of Computational Science}
  \bibinfo{volume}{65}, \bibinfo{pages}{101883}.
\newblock \DOIprefix\doi{10.1016/j.jocs.2022.101883}.
\bibitem[{Verner(2010)}]{verner2010numerically}
\bibinfo{author}{Verner, J.H.}, \bibinfo{year}{2010}.
\newblock \bibinfo{title}{Numerically optimal {Runge}--{Kutta} pairs with
  interpolants}.
\newblock \bibinfo{journal}{Numerical Algorithms} \bibinfo{volume}{53},
  \bibinfo{pages}{383--396}.
\newblock \DOIprefix\doi{10.1007/s11075-009-9290-3}.
\bibitem[{Wittig et~al.(2015)Wittig, Di~Lizia, Armellin, Makino,
  Bernelli-Zazzera and Berz}]{wittig2015propagation}
\bibinfo{author}{Wittig, A.}, \bibinfo{author}{Di~Lizia, P.},
  \bibinfo{author}{Armellin, R.}, \bibinfo{author}{Makino, K.},
  \bibinfo{author}{Bernelli-Zazzera, F.}, \bibinfo{author}{Berz, M.},
  \bibinfo{year}{2015}.
\newblock \bibinfo{title}{Propagation of large uncertainty sets in orbital
  dynamics by automatic domain splitting}.
\newblock \bibinfo{journal}{Celestial Mechanics and Dynamical Astronomy}
  \bibinfo{volume}{122}, \bibinfo{pages}{239--261}.
\newblock \DOIprefix\doi{10.1007/s10569-015-9618-3}.

\end{thebibliography}

\end{document}